\documentclass[epic,11pt]{article}
\usepackage{amsfonts,amsmath,amscd,amssymb,amsthm}
\usepackage{fullpage}
\usepackage[alphabetic,initials]{amsrefs}
\usepackage[all]{xy}
\usepackage{hyperref} 
\usepackage{mathrsfs}
\input xy%
\usepackage{fullpage}
\usepackage[pdftex]{pict2e}
\def\0{{\bar 0}}
\def\1{{\bar 1}}

\def\Z{{\mathbb Z}}
\def\X{{\mathbb X}}
\def\Y{{\mathbb Y}}
\def\F{{\mathbb F}}

\def\cpa{{\operatorname{c.empty}}} 
\def\cda{{\operatorname{c.full}}} 
\def\rpa{{\operatorname{r.empty}}} 
\def\rda{{\operatorname{r.full}}}

\def\W{\mathtt{W}}

\def\ua{{\operatorname{a  }}}
\def\ur{{\operatorname{r  }}}
\def\ud{{\operatorname{d}}}

\def\str{{\operatorname{red}}}
 
\def\trans{{\operatorname{transpose}}} 

\def\iso{{\operatorname{iso}}}

\newcommand{\ttk}{\mathtt{k}}
\newcommand{\ttr}{\mathtt{r}}
\newcommand{\ttw}{\mathtt{w}}
\newcommand{\ttc}{\mathtt{c}}

\newcommand{\ttd}{\mathtt{d}}

\newcommand{\lst}{{\stackrel{\leftarrow}{\gl}}}

\newcommand{\bA}{{\bf a}}

\newcommand{\diag}{diag}

\newcommand{\Lgh}{\widehat{L}(\stackrel{{\rm _o}}{\fg})}

\newcommand{\itema}{\item[{{\rm$($a$)$}}]}
\newcommand{\itemb}{\item[{{\rm$($b$)$}}]}
\newcommand{\itemc}{\item[{{\rm$($c$)$}}]}
\newcommand{\itemd}{\item[{{\rm$($d$)$}}]}
\newcommand{\iteme}{\item[{{\rm$($e$)$}}]}

\newcommand{\itemo}{\item[{}]}
\newcommand{\noi}{\noindent}

\newcommand{\ga}{\alpha}
\newcommand{\gb}{\beta}

\newcommand{\Gc}{\Gamma}
\newcommand{\Gl}{\Lambda}
\newcommand{\Gd}{\Delta}

\newcommand{\gd}{\delta}

\newcommand{\go}{\omega}

\newcommand{\gt}{\tau}

\newcommand{\gl}{\lambda}

\newcommand{\gr}{\rho}
\newcommand{\gk}{\kappa}
\newcommand{\gep}{\epsilon}

\newcommand{\A}{\mathbb A}

\newcommand{\mr}{{-_\ttr}}
\newcommand{\pr}{{+_\ttr}}

\newcommand{\omc}{{-_\ttc}}
\newcommand{\opc}{{+_\ttc}}

\def\Im{{\operatorname{Im}\;}}

\newcommand{\fg}{\mathfrak{g}}\newcommand{\fgl}{\mathfrak{gl}}

\newcommand{\fh}{\mathfrak{h}}

\newcommand{\fB}{\mathfrak{B}}

\newcommand{\fT}{\mathfrak T}

\newcommand{\ff}{\footnote}
\newfont{\eufm}{eufm10 scaled\magstep1}

 \newcommand{\ti}{\times}

\newcommand{\bcu}{\bigcup}

\newcommand{\bsk}{\backslash}

\newcommand{\cO}{\mathcal{O}}

\newcommand{\cG}{\mathcal{G}}

\newcommand{\cS}{\mathcal{S}}

\newcommand{\cX}{\mathcal{X}}

\newcommand{\cW}{\mathfrak{W}}

\newcommand{\ey}{\end{eqnarray}}
\newcommand{\by}{\begin{eqnarray}}
\newcommand{\nn}{\nonumber}

\newcommand{\bco}{\begin{conjecture}}
\newcommand{\ba}{\begin{alg}}
\newcommand{\ea}{\end{alg}}
\newcommand{\eco}{\end{conjecture}}
\newcommand{\bpf}{\begin{proof}}
\newcommand{\epf}{\end{proof}}
\newcommand{\bt}{\begin{theorem}}

\newcommand{\et}{\end{theorem}}

\newcommand{\br}{\begin{rem}}
\newcommand{\er}{\end{rem}}
\newcommand{\brs}{\begin{rems}}
\newcommand{\ers}{\end{rems}}
\newcommand{\bi}{\begin{itemize}}
\newcommand{\ei}{\end{itemize}}
\newcommand{\bl}{\begin{lemma}}
\newcommand{\bsul}{\begin{sublemma}}
\newcommand{\esul}{\end{sublemma}}
\newcommand{\bp}{\begin{proposition}}
\newcommand{\be}{\begin{equation}}
\newcommand{\bc}{\begin{corollary}}
\newcommand{\bexs}{\begin{examples}}
\newcommand{\eexs}{\end{examples}}
\newcommand{\bexa}{\begin{example}}
\newcommand{\eexa}{\end{example}}
\newcommand{\bex}{\begin{exercise}}
\newcommand{\eex}{\end{exercise}}
\newcommand{\btab}{\begin{tab}}
\newcommand{\etab}{\end{tab}}
\newcommand{\el}{\end{lemma}}
\newcommand{\ep}{\end{proposition}}
\newcommand{\ee}{\end{equation}}
\newcommand{\ec}{\end{corollary}}
\newcommand{\Bc}{\begin{center}}
\newcommand{\Ec}{\end{center}}
\newcommand{\bh}{\begin{hyp}}
\newcommand{\eh}{\end{hyp}}
\newcommand{\bhs}{\begin{hyps}}
\newcommand{\ehs}{\end{hyps}}
\newcommand{\bd}{\begin{dfn}}
\newcommand{\ed}{\end{dfn}}

\newcommand{\bn}{\begin{notn}}
\newcommand{\en}{\end{notn}}

\newcommand{\fgo}{\stackrel{{\rm o}}{\fg}}

{}
{}

\newcommand{\fho}{\stackrel{\rm o}{\fh}}
\newcommand{\fhs}{\stackrel{{\rm _o}}{\fh}^*}

\usepackage[aligntableaux=center]{ytableau}
\usepackage{tabularray}

\begin{document}
\title{Table of Contents}


\newtheorem*{bend}{Dangerous Bend}

\newtheorem{thm}{Theorem}[section]
\newtheorem{hyp}[thm]{Hypothesis}
 \newtheorem{hyps}[thm]{Hypotheses}
\newtheorem{notn}[thm]{Notation}

  \newtheorem{rems}[thm]{Remarks}

\newtheorem{conjecture}[thm]{Conjecture}
\newtheorem{theorem}[thm]{Theorem}
\newtheorem{theorem a}[thm]{Theorem A}
\newtheorem{example}[thm]{Example}
\newtheorem{examples}[thm]{Examples}
\newtheorem{corollary}[thm]{Corollary}
\newtheorem{rem}[thm]{Remark}
\newtheorem{lemma}[thm]{Lemma}
\newtheorem{sublemma}[thm]{Sublemma}
\newtheorem{cor}[thm]{Corollary}
\newtheorem{proposition}[thm]{Proposition}
\newtheorem{exs}[thm]{Examples}
\newtheorem{ex}[thm]{Example}
\newtheorem{exercise}[thm]{Exercise}
\numberwithin{equation}{section}%
\setcounter{part}{0}
\newcommand{\drar}{\rightarrow}
\newcommand{\lra}{\longrightarrow}
\newcommand{\rra}{\longleftarrow}
\newcommand{\dra}{\Rightarrow}
\newcommand{\dla}{\Leftarrow}
\newcommand{\rl}{\longleftrightarrow}

\newtheorem{Thm}{Main Theorem}


\newtheorem*{thm*}{Theorem}
\newtheorem{lem}[thm]{Lemma}
\newtheorem*{lem*}{Lemma}
\newtheorem*{prop*}{Proposition}
\newtheorem*{cor*}{Corollary}
\newtheorem{dfn}[thm]{Definition}
\newtheorem*{defn*}{Definition}
\newtheorem{notadefn}[thm]{Notation and Definition}
\newtheorem*{notadefn*}{Notation and Definition}
\newtheorem{nota}[thm]{Notation}
\newtheorem*{nota*}{Notation}
\newtheorem{note}[thm]{Remark}
\newtheorem*{note*}{Remark}
\newtheorem*{notes*}{Remarks}
\newtheorem{hypo}[thm]{Hypothesis}
\newtheorem*{ex*}{Example}
\newtheorem{prob}[thm]{Problems}
\newtheorem{conj}[thm]{Conjecture}

\title{ Young diagrams, deformed Calogero-Moser systems and Cayley graphs}
\author{Ian M. Musson
\\Department of Mathematical Sciences\\
University of Wisconsin-Milwaukee\\ email: {\tt
musson@uwm.edu}}
\maketitle
\begin{abstract} Let $\ttk$ be an algebraically closed field of characteristic zero and $n, m$ coprime positive integers.  
Let $\fgo$ be the Lie superalgebra $\fgl(n|m)$ with root system $\Gd$. Using $\Gd$,   Sergeev and Veselov, \cite{SV2} introduced an action of the Weyl groupoid 
$\cW$, in connection with their  study of the 
the Grothendieck group of finite dimensinonal graded $\fg$-modules.  
We denote the subgroupoid of  $\cW$ with morphisms corresponding to isotropic roots by
$\mathfrak T_{iso}$.  

Later, \cite{SV101}  the same authors defined an action of $\cW$ on $X=\ttk^{n|m}$ such that the invariant algebra $\cO(X)^\cW$  
is isomorphic  to the algebra of quantum integrals for the deformed Calogero-Moser  system 
introduced in \cite{SV1}. This completely integrable system depends on a non-zero parameter $\gk$. When $\gk=-m/n$ we study a certain infinite $\mathfrak  T_{iso}$-orbit    {\bf O}  for this action. 
The 
Cayley graph  for this orbit is isomorphic to  the Cayley graphs for two other actions of 
$\mathfrak T_{iso}$  which were studied in  \cite{M23}.   

\end{abstract}



\section{Introduction} \label{akmls} 

For a set $\X$, the symmetric groupoid $\cS(\X)$  on  $\X$ has base consisting of all subsets of $\X$ and morphisms all bijections between subsets.
The groupoid $\cG$ {\it  acts on} $\X$ if there is a functor $F:\cG\lra\cS(\X)$.  If $X$ is an affine algebraic variety, we can study the invariant algebra $\cO(X)^\cG$,  consisting of all polynomial functions on $X$ that are constant on $\cG$-orbits, 
\cite{SV101}, \cite{M22}.   
Let  $\cW$ be the Weyl groupoid in type  $A(n, m)$ from  \cite{SV2}. 
We denote the subgroupoid of  $\cW$ with morphisms corresponding to isotropic roots by
$\mathfrak T_{iso}$.  
\\ \\
In  \cite{SV1} 
 Sergeev and Veselov studied  deformed quantum Calogero-Moser
 (CM)   systems related to generalized root systems (in the sense of Serganova \cite{S}).  
For the classical types $A(n, m)$ and $BC(n, m)$, the admissible deformations depend on a single parameter and are completely 
integrable.  
Later \cite{SV101} in type A,  they described the algebra of quantum integrals for the deformed CM system as the the invariant algebra $\cO(X)^\cW$  for an action of $\cW$ on $X=\ttk^{n|m}$.  
\\ \\
We consider a certain infinite $\mathfrak  T_{iso}$-orbit    {\bf O}  which appears in \cite{SV101} Equation (14).  Up to some inherent symmetries, this is the unique infinite orbit for the parameters we consider, see Conjecture \ref{t2}.
The 
Cayley graph  for this orbit is isomorphic to  the Cayley graphs for two other actions of 
$\mathfrak T_{iso}$.   The first arises from the set of Young diagrams $X$ contained in a rectangle with $n$ rows and $m$ columns. 
We can define an equivalence relation on $X\ti \Z$ such that   $\mathfrak  T_{iso}$ acts on the  set of equivalence classes, see Subsection \ref{itr3}.  Secondly  $\mathfrak  T_{iso}$ acts transitively, by odd reflections on  a set of  Borel subalgebras of  the affinization $\Lgh$ of  $\fgo$, see \cite{M23} for details.
\subsection{ Groupoids and Cayley graphs }\label{itr1} 
Throughout  $[k]$ denotes  the set of the first $k$ positive integers and iff means if and only if.                                                                                                                                                                                                                                                                                                                                             
Suppose $F$ is a functor $F:\cG\lra\cS(\X)$.
The {\it Cayley graph} $\Gc(F)$ of $F$ is the 
directed graph with vertex set $\X$ and an edge $x\lra F(s)x$ whenever $s$ is a morphism in $\cG$ and $F(s)x$ is defined.  
We assume $\cG=\mathfrak  T_{iso}$. 
 If $x\in \X$, we refer to the (vertices in the) connected component of the Cayley graph as the {\it orbit} of $x$.  
We say that $F$ acts {\it transitively} if $\X$  is a single $\cG$-orbit.
\\ \\
If $\go:\X\lra \Y$ is a bijection of sets, there is an induced isomorphism $\cS(\X)\lra \cS(\Y)$  also denoted $\go$.  
   Now suppose we have a commutative diagram of groupoids and functors
\be \label{bby}
\xymatrix@C=2pc@R=1pc{
\cG\ar@{<->}_= [dd]&&
\cS(\X)\ar@{<-}_F[ll]: \ar@{->}^{\go}[dd]&\\ \\
\cG \ar@{-}[rr]_B &&
\cS(\Y)\ar@{<-}[ll] &}
\ee
Thus $\go$ induces an isomorphism of Cayley graphs 
$$\Gc(F) \lra \Gc(B).$$ If $F$ and $B$ are transitive in \eqref{bby}, we say that $\go$ is a $\cG$-{\it equivariant bijection from 
$\X$ to $\Y$}. 
The aim of this paper is to construct various commutative diagrams as in \eqref{bby}.

\subsection{The Lie superalgebra $\fgl(n|m)$ and the Weyl groupoid
}\label{Lsg}
Let $\fgo$ be the Lie superalgebras $\fgl(n|m)$ with  root system $\Gd$. The Cartan subalgebra $\fho$ of diagonal matrices in $\fgo$ is identified with $\ttk^{n|m}$.  
For $i\in [m+n], $  define $\gep_i\in \fhs$ so that $\gep_i(h) = h_i$ where  $h=\diag(h_1,\ldots,h_{m+n})\in \;\fho$.  For $i\in [m]$, set $i'=i+n$ and  $\gd_{i} =\gep_{i'}$. 
The bilinear form $(\;,\;)$ is defined on $\fh^*$ is defined on the basis 
$\epsilon_{1},\ldots,   \epsilon_{n}, \delta_{1}, \ldots,\delta_{m}  $ by
\be \label{edform}(\epsilon_i,\epsilon_j) = \delta_{i,j} = - (\delta_i,\delta_j)\nn\ee 
and $(\epsilon_{i},\gd_{j})=0$ for all relevant indices $i,j$.
We write  $\Gl \in \ttk^{n|m}$ in the form
\be \label{ezr} \Gl = 
(a_1, \ldots ,a_n|b_1,\ldots,b_{m}) =  \sum_{i=1}^n  a_i
\gep_i - \sum_{j=1}^m  b_j \gd_j .\ee
If $\Gl$ is as above we have 
\be \label{exy} (\Gl,\gep_i- \gd_j)=  a_i-b_j. \ee 
Let $\mathfrak T_{iso}$ be the groupoid with base 
$\Gd_{iso}$, the set of all the isotropic roots in $\Gd$  
and  non-identity  morphisms  $\gr_{\alpha}:\alpha \rightarrow -\ga$, $\ga\in\Gd_{iso}$. The Weyl group $W$ acts on $\mathfrak T_{iso}$ in a natural way: $\alpha \rightarrow w(\alpha),\,
\gr_{\alpha} \rightarrow \gr_{w(\alpha)}$.  The {\it Weyl groupoid} $
\cW$ is defined
by 
\be \label{fgh}  \mathfrak{W} = W \coprod W \ltimes \mathfrak T_{iso},\ee     the disjoint union of the group $W$ considered as a groupoid with a single point base $[W]$ and the semidirect product $W \ltimes \mathfrak T_{iso}$. 
If $\ga\in \Gd^+_1$ we call $\gr_\ga$ a {\it positive morphism} and likewise for the morphisms $t_\ga =  F(\gr_\ga)$, $\F(\gr_\ga) =T_\ga$ 
and $\gt_\ga =  SV(\gr_\ga)$ defined in \eqref {BPG}, Theorem \ref{shg} and Subsection \ref {aog11}.

\subsection{ Partitions} \label{itr2} 
A {\it partition} $\gl$ into $n$ {\it parts} is a  descending sequence 
\be \label{pdef}  
\gl=(\gl_1, \gl_2, \ldots, \gl_n)\ee of non-negative integers.  
We call $\gl_i$ the $i^{th}$ {\it part} of  $\gl$.  To minimize notation, a
  partition is identified with the Young diagram it defines. 
Given positive integers $m, n$, 
let $X$ be the set of Young diagrams that fit inside (the bottom left corner of) 
a  rectangle {\bf R}  with $n$ rows and $m$ columns.  Boxes in 
{\bf R} correspond to root spaces in  the 
upper triangular block $\fg^+_1$ of $\fgl(n|m)$.  Accordingly the rows of   
{\bf R} are indexed from the top down by $\gep_1,\ldots,\gep_n$  and columns of {\bf R} are indexed  from  left to right by  $\gd_1,\ldots,\gd_m$.  Thus the box in row $\gep_i$ and column $\gd_j$
corresponds to the root space $\fg^{\gep_i-\gd_j}$ and we call this box $\fB(\gep_i-\gd_j)$, or often just $\gep_i-\gd_j$.  The set of roots of 
$\fg^+_1$ is $\Gd_1^+ = \{\gep_i-\gd_j|i\in [n], j \in [m]\}.$
\\ \\
We also write 
\be \label{prf}\gl=(m^{a_m},  \ldots, i^{a_i}, \ldots)\ee
 to indicate that the part $i$ occurs with multiplicity ${a_i}$ in $\gl$.   Set $|\gl| = \sum_{i=1}^n \gl_i$. 
For a partition $\gl$,  the {\it dual partition} $\gl'$  has $j^{th}$ part 
$\gl'_j = |\{i|\gl_i \ge j\}|$. The Young diagram corresponding to $\gl$ has $\gl_i$ boxes in row $\gep_{n+1-i}$ with the first box on the left edge of {\bf R}. Our convention differs from the usual labelling of rows in a Young diagram, however we will also write matrices in {\bf R}, and it agrees with the usual way of counting rows in a matrix. To    
resolve this minor issue we  define the vector $\lst$ in \eqref{cd}, reversing the entries in $\gl$.
Denote the empty Young diagram $(0^n)$ by $\emptyset$.
\\ \\
In what follows we view a Young diagram $\gl$ as a set of boxes in {\bf R}.
If $\ga\in\Gd_1^+$,we say the box  $\fB=\fB( \ga)$ is an  {\it outer corner} of $\gl$ if $
  \fB \notin \gl$ 
and the addition of $ \fB$ to $\gl$ creates a new Young diagram  $t_{ \ga}(\gl) $
Similarly box $ \fB$ is an  {\it inner corner} of $\gl$ if $
 \fB \in \gl$ 
and the removal of $ \fB$ from $\gl$ creates a new Young diagram 
 $t_{- \ga}(\gl)$.  
Let 
$$X_\ga = \{\gl\in X| \fB(\ga) \mbox{ is an outer corner  of } \gl\}$$
and  
$$X_{-\ga} = \{\gl\in X| \fB(\ga) \mbox{ is an inner corner  of } \gl\}.$$  
Then   $t_{ \pm\ga} $ define mutually inverse bijections  $ X_{\ga} \rl X_{-\ga} $ and 
we have  a functor $F:\fT_{\iso}\lra \cS(X)$ given by 
\be \label{BPG}  F({\alpha})=X_\ga \mbox{  and } F(\gr_{\alpha}) =t_\ga.\ee
\subsection{ The equivalence relation on $X\ti \Z$} 
\label{erx}
We say that $\gl \in X$ is {\it row empty} if the highest row of $\gl$ is empty 
 and {\it row full} if the lowest row of $\gl$ contains $m$  boxes. The terms {\it column empty} 
 and {\it column full} are defined similarly. More formally \be\label{yyt}X(\rda) = \{\gl\in X|\gl_1=m\}, \quad X(\rpa) = \{\gl\in X|\gl_n=0\}.\ee
\be\label{yys}X(\cda) = \{\gl\in X|\gl'_1=n\}, \quad X(\cpa) = \{\gl\in X|\gl'_m=0\}.\ee
If $\gl \in X(\rda)$, let
 $\gl^\mr$ be the Young diagram with row $\gep_{1}$ empty and such that 
for $i\in [n-1] $ row $\gep_{i+1}$ of 
 $\gl^\mr$   equals row $\gep_{i}$ of $\gl$.   There is a bijection 
$X(\rda) \lra X(\rpa)$ given by  $\gl \lra \gl^\mr$ whose inverse we denote by  
$ \mu\lra \mu^\pr.$
Similarly, if $\gl \in X(\cda)$, 
let   $\gl^\omc$ be the Young diagram obtained from $\gl$  by deleting the first column.     We have  inverse bijections
$X(\cda) \rl X(\cpa)$ given by  $\gl \lra \gl^\omc$ and 
$\mu^\opc \rra \mu.$
\noi Let $\sim$ be the smallest equivalence relation on 
$X\ti \Z$ such that 
\be \label{1kr}(\gl,k)\sim ({\gl}^\mr,k+m) \mbox{ if } \gl\in X(\rda)\ee and 
\be \label{2kr}(\gl,k)\sim ({\gl}^\omc,k+n) \mbox{ if } \gl\in X(\cda).\ee
Denote the equivalence class of $(\gl,k)\in X\ti Z$ under $\sim$ by $[\gl,k]$. 
The action of $\mathfrak T_{iso}$  can be extended to the set of equivalence classes  $[X \ti \Z]$, see Subsection \ref{itr3}. 
\subsection{ Main result} \label{itr6}

We reformulate the definition of  action of the   $\mathfrak{W}$  on $\ttk^{n|m}$ from  \cite{SV101} to obtain 
 a functor 
\be \label{e11}SV:\mathfrak  T_{iso}\lra \mathfrak \cS(\Z^{n|m}).\ee
For $\ga=\gep_i- \gd_j \in \Gd_\iso^+$, 
 define $v_\ga = n\gep_i -m \gd_j$ and 
\be \label{e99}\Pi_{\ga} = \{\Gl\in \ttk^{n|m}| (\Gl, \ga)=0\},\quad \Pi_{-\ga} = \{\Gl\in \ttk^{n|m}| (\Gl,\ga)=n-m\}.\ee
Next define $\gt_{\pm\ga} :\Pi_{\pm\ga} \lra \Pi_{\mp\ga} $ by 
\be \label{e9}\gt_{\pm\ga} (\Gl) = \Gl \pm v_\ga.\ee 
The 
functor  $SV$ from \eqref{e11} is given  by 
$SV(\pm \ga) =\Pi_{\pm\ga}$, and $SV(\gr_{\pm\ga}) =\gt_{\pm\ga}$. 
Using the permutation   action of $W$, we obtain  an action of $\cW$ on $\Z^{n|m}$ (and also on $\ttk ^{n|m}$).
More details are given in Subsection \ref{aog11}. 
Let {\bf O}
be the  $\mathfrak  T_{iso}$-orbit  
of 
$\Gl_0 \in \Z^{n|m}$  from   \eqref{e7}.
\bt \label{iir}  Assume  $m, n$ are coprime.  
The map $\hat x$  from Theorem \ref{rjt}  is a $\mathfrak T_{iso}$-equivariant bijection $\hat x:[X \ti \Z] \lra{\bf O}$ such that $\hat x[\emptyset, 0] =   \Gl_0$. 
\et \noi  
This paper is organized as follows. Section \ref{fc} contains preliminary results on Young diagrams  and related combinatorics.  
Orbits for the action of $\mathfrak  T_{iso}$  on $\Z^{n|m}$ 
are studied in 
Section \ref{affc}.   
We end with some concluding remarks in  
Section \ref{cio}.  
In Examples 
\ref{E2} we illustrate many of our results in the  case
where $(n,m)=(2,3)$.
\section{Background Results.
}\label{fc}
\subsection{Rotation operators}\label{roo}

Define the permutation  $\nu$ of $\{1, \ldots, n\}$ by 
$\nu(k) = k-1 \mbox{ mod } n$. 
We extend $\nu$ to  an operator on $\ttk^{n|m}$  by setting
\be \label{e67l}\nu(a_1, \ldots, a_n|b_1, \ldots, b_m) = (a_{\nu(1)}, \ldots, a_{\nu(n)}|b_1, \ldots, b_m) \ee
Similarly define the permutation 
$\eta $ of $\{1', \ldots, m'\}$ by 
$\eta(k') = (k+1)' \mbox{ mod } m$ and set  
\be \label{e67l}\eta(a_1, \ldots, a_n|b_1, \ldots, b_m) = (a_{ 1}, \ldots, a_{n}|b_{\eta(1)}, \ldots, b_{\eta(m)}) \ee
Now \eqref{e67l} implies 
\be \label{e67m}\nu(\gep_i) = \gep_{i+1}  \mbox{ for } i\in [n], \quad \eta(\gd_{j+1}) = \gd_{j}  \mbox{ for } j\in [m].\ee
where the subscripts are read mod $n$ on the left and mod $m$ on the right, and 
$$\nu(\gd_j) = \gd_j ,\quad \eta(\gep_i)=\gep_{i} 
\mbox{ for } j\in [m], \quad 
 i\in [n].$$

\bl
If  $\gb= { \eta^{i} \nu^{j}}\ga$ and $\Gl\in \Pi_{\pm\gb}$ 
then 
${ \eta^{-i} \nu^{-j}} \Gl\in \Pi_{\pm\ga}$ 
and 
\be \label{bc1} 
\gt_{\pm\ga}
{ \eta^{-i} \nu^{-j}}\Gl
=
{ \eta^{-i} \nu^{-j}}\gt_{\pm\gb}\Gl
.\ee
\el 
\bpf Let $\ga = \gep_{r} - \gd_{s}$ and $\gb = \gep_{r+j} - \gd_{s-i}$.  For positive morphisms the result follows from the display below in which the first row gives the positions of the entries in the lower rows.  We assume that $ {r+j} \le n$ and ${s-i} \ge 1$.  Otherwise since the indices on $\gep_{r+j} $ and  $\gd_{s-i}$  are read mod $n$ and $m$ respectively, the order of  the entries in the rows is changed, but the proof is essentially the same.
\by &&\quad \quad \quad  r  \quad \quad r+j \quad  \quad  \quad   \quad  (s-i)'  \quad \quad s'  \quad \quad  \quad   \nn
\\ \Gl &=&( \;\ldots, \; \ldots,  \; \ldots, \; a, \; \ldots \;|  \;\;\ldots,  \;a,   \;\;\ldots,  \; \;\ldots,  \;\ldots  \; )\nn\\
\eta^{-i} \nu^{-j} \Gl &=& ( \;\ldots,  \; a, \; \ldots, \; \ldots, \; \ldots \;| \;\ldots,  \; \ldots, \; \ldots,   \;\;a,  \; \; \;\ldots   \; \;\; )\nn\\
\gt_{\gb}\Gl &=& (\ldots,  \ldots, \ldots, a+ n, \ldots|\ldots, a+ m, \ldots, \ldots, \ldots\;)\nn\\
\gt_{\ga}{ \eta^{-i} \nu^{-j}}\Gl
=
{ \eta^{-i} \nu^{-j}}\gt_{\gb}\Gl&=& (\ldots,  a+ n, \ldots, \ldots, \ldots| \;\ldots,  \ldots, \ldots, a+ m, \ldots\;).\nn
\ey
\epf
\subsection{Young diagrams,  partitions and corners}\label{ypc}
The next two Lemmas are proved in \cite{M23}.
\bl \label{rjs} 
If $\ga = \gep_i-\gd_{j}$, the following are equivalent
\bi \itema $ \fB(\ga)$ is an inner corner  of $\gl$, that is $\gl\in X_{-\ga}$.
\itemb $\gl_{n+1-i} = j$ 
and $\gl'_{j} = n+1-i$ .
\ei
If these conditions hold then 
\be \label{qef}  t_{-\ga}(\gl) = (\gl_1,\ldots , \gl_{n-i}, \gl_{n+1-i}-1,\gl_{n+2-i},\ldots, \gl_n).\ee
\el
\bl \label{rjx} 
If $\ga = \gep_i-\gd_{j}$, the following are equivalent
\bi \itema $ \fB(\ga)$ is an outer corner  of $\gl$, that is $\gl\in X_{\ga}$.
\itemb $\gl_{n+1-i} = j-1$ and $\gl'_{j} = n-i$.
\ei
If these conditions hold then 
\be \label{qhf} t_{\ga}(\gl) = (\gl_1,\ldots , \gl_{n-i}, \gl_{n+1-i}+1,\gl_{n+2-i},\ldots, \gl_n).\ee
\el \noi 
The partition 
$ \boldsymbol \gl  = (m,1^{n-1})$ plays an important role.
 The set of {\it reduced diagrams} is 
$$X^{\str}= \{\gl \subseteq {\bf R} |  \gl_n = \gl'_m =0\}.$$ \ff{ 
If $\gl\in X$ let $\bar \gl \in X$ be the complement of the boxes obtained from $\gl$ by  rotating $X$ through angle $\pi$.  The involution $\gl \rl \bar \gl$ interchanges 
$X^{\str}$ and 
$\{\gl\in X| \boldsymbol \gl  \subseteq \gl\}$.  This remark can be used to explain the existence of many parallel results for inner and outer (psuedo-corners).}
If ${ \boldsymbol \gl} \subseteq \gl$,
we call $\gep_n-\gd_1$   a   {\it outer pseudo-corner} of  $\gl$, and if 
$\gl \in X^{\str}$ we call $ \gep_1-\gd_m  $   a {\it  inner  pseudo-corner} of  $\gl$.  These are not outer  or inner corners, but behave in some ways as if they were.  Pseudo-corners 
play an important role in  Section \ref{affc}, see 
Lemmas \ref{c1} and \ref{c2}.  
It is easy to see the following.
\bl \label{rey} \bi \itemo
\itema  If 
$\ga =\gep_n-\gd_1$ 
and ${ \boldsymbol \gl} \subseteq  \gl$,
 then $\ga$ is not an outer corner of $\gl$, but 
$\nu\ga = \gep_1-\gd_1$ is an outer corner of 
${\gl}^\mr$ and $\eta^{}\ga=\gep_n 
-\gd_m$ is  an outer corner of ${\gl}^\omc.$
\itemb  
If $\ga =\gep_1-\gd_m$ and $\gl \in X^{\str}$, then $\ga$ is not an inner corner of $\gl$, but $\nu^{-1}\ga=\gep_n 
-\gd_m$ is  an inner corner of ${\gl}^\pr$ and $\eta^{-1}\ga=\gep_1 
-\gd_1$ is  an inner corner of ${\gl}^\opc.$
\ei
\el  \noi

 \subsection{Words}\label{1bs}

 \subsubsection{Young diagrams  and words}\label{bsa2}
To a Young diagram $\gl$ in $X$, we associate a word $\ttw=\ttw(\gl)$ of length $m+n$
in the alphabet $\{\ttr, \ttd\}$.  The upper border of $\gl$ in $X$ is a connected sequence of $m+n$ line segments of unit length starting at the 
top left corner and ending at the bottom right corner of
 {\bf R}.  Each line segment carries the label
$\ur$ or $\ud$, according to whether it moves to the right or down.  Then $\ttw$ is the sequence of labels.  Let $\W$ be the set of words with    
$m$ $\ur$'s and $n$ $ \ud$'s.  There is a bijection  $X\lra \W$ given by
$\gl \lra\ttw(\gl)$.
\subsubsection{ The equivalence relation} \label{tgs}
Using words it is easy to describe the equivalence class of $(w_0, k_0)$ under $\sim$,  where $w_0 = \ttw(\gl)$.  We recall some results from  \cite{M23}. 
Inductively define words $w_i$ ($i>0$) using the following rules 
\bi
\itema If  $w_i=\ttd x$, then $w_{i+1} = x\ttd $ and $k_{i+1} = k_i -m$.
\itemb If  $w_i=\ttr x$, then $w_{i+1} = x\ttr $ and $k_{i+1} = k_i +n$.
\ei 
Define $\gl_i \in X$ by   $w_{i} =\ttw (\gl_i)$.
 \bl \label{iod}
 \begin{itemize}\itemo
\itema If the first letter in $w_i$ is $\ttr$, then $\gl_{i+1}=\gl_{i}^\omc$.
\itemb If the first letter in $w_i$ is $\ttd$, then $\gl_{i+1}=\gl_{i}^\pr$.
\ei
\noi 
 \el
\noi
Now each word contains $m$ symbols $\ttr$ and $n$ symbols $\ttd$.  Therefore $w_{m+n} = w_0.$
To find a closed formula for $k_i$ set 
\[ \begin{array}
  {rcl}
 c_j =1, d_j =0 & \mbox{if} & \mbox{the } j^{th}  \mbox{ letter of } \ttw \mbox{ is } \ttr\\
 c_j =0, d_j =1 & \mbox{if} & \mbox{the } j^{th}  \mbox{ letter of } \ttw \mbox{ is } \ttd\end{array},\] then define
\be \label{ips}  k_{i} = k_0 +n\sum_{j=1}^{i}c_j -m\sum_{j=1}^{i}d_j. 
\ee

\bl \label{vgs} \itemo\bi \itema 
The set $E:=\{(w_{i},k_{i})| 0 \le  i < m+n\}$, is the  equivalence class of $(w_0,k_0)$ in $[\W\ti\Z]$. 
\itemb   The  equivalence class of $(\gl_0,k_0)$ in $[X\ti\Z]$ is $\{(\gl_{i},k_{i})| 0 \le  i < m+n\}$. 
\ei
\el

\subsection{ Extending  the action of $\mathfrak T_{iso}$} 
\label{itr3}
We study the interaction between adding (or deleting) rows and  corners.  We can see what to expect by looking at small examples.
\bexa  {\rm 
Below a full row consists of $m=3$ boxes, $n\ge 2$, $\gb = \gep_{n-1} -\gd_2$ and $\ga = \nu\gb= \gep_{n} -\gd_2$.  
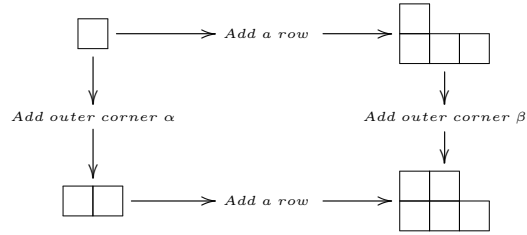
\tiny
\begin{equation*}
\xymatrix@C=1pc@R=1pc{ &
\ydiagram{1}
\ar@{->}[d]&  \scriptscriptstyle{Add\;  a  \; row\;}\ar@{<-}[l] \ar@{->}[r]& \ydiagram{0+1,3}
& 
&& &
\\&\scriptscriptstyle{Add\;  outer \;  corner\;} \ga \ar@{->}[d]&&\scriptscriptstyle{Add\;  outer \;  corner\;} \gb \ar@{<-}[u]& &
\\
&\ydiagram{2}\ar@{->}[r] & \scriptscriptstyle{Add\;  a  \; row\;}\ar@{->}[r] &
\ydiagram{0+2,3}\ar@{<-}[u]&& &&}
\end{equation*}
\large
\normalsize
}\eexa 
\noi
More generally this shows part (a) of the following result.  
\bl \label{rhy} 
\bi \itemo \itema
If
$\gl\in  X(\rpa)\cap X_{\nu\gb}$, then
$\gl^{+_\ttr}\in X_{\gb}$, ${t_{\nu\gb}(\gl)}\in X(\rpa)$ and  
\be \label{kxt} {t_{\gb}(\gl}^{+_\ttr}) = t_{\nu\gb}({\gl})^{+_\ttr}.\ee
\itemb If
$\gl\in  X(\cpa)\cap X_{\eta\gb}$, then  
$\gl^{+_\ttc}\in X_{\gb}$ , ${t_{\eta\gb}(\gl)}\in X(\cpa)$ and  
\be \label{kzt} {t_{\gb}(\gl}^{+_\ttc}) = t_{\eta\gb}({\gl})^{+_\ttc}.\ee
\ei\el
\bpf For a proof without diagrams see \cite{M23}.
\epf
 \noi 
Assume  $\mathfrak  T_{iso}$  acts on $X$ as in \eqref{BPG}  and  that $n,m$ are coprime.  We  define  an action of $\mathfrak T_{iso}$ on the set of equivalence classes $\cX=[X\ti Z]$ under the equivalence relation generated by \eqref{1kr} 
and \eqref{2kr}. This was done in greater generality in \cite{M23}.   Here we can  make certain simplifications since we do  not consider Borel subalgebras.
First set 
\be \label{sy1}
F_{in+jm } (\ga)= (X_{ \eta^{i}\nu^{j}\ga},in+jm )
\ee and let $\cX_\ga$ 
be the set of equivalence classes of elements in $\bcu_{i, j \in \Z} F_{in+jm } (\ga).$ 
\noi 
\bt \label{shg}  
 There is a functor $\F:\mathfrak  T_{iso}\lra \cS[X\ti \Z]$ given by $\F(\ga) = \cX_\ga $ 
and
\be \label{dyx}
\F(\gr_\ga)[\gl, in+jm] =
 [t_{\eta^{i}\nu^{j}\ga}(\gl),in+jm]
\mbox{ for  } (\gl, in+jm)\in
F_{in+jm } (\ga)
.\ee 
\et \noi 
Set  $\F(\gr_\ga) =T_\ga$. 
Note that $ (\gl^{+_\ttr}, in+jm) \sim  (\gl, in+(j+1)m)$.  Thus to show \eqref{dyx} is well-defined we need to show.
\bl If $\gl$ and $ \gb=\eta^{i}\nu^{j}\ga$  satisfy the hypothesis of Lemma \ref{rhy} (a), then 
$$ T_\ga(\gl^{+_\ttr}, in+jm) \sim  T_\ga(\gl, in+(j+1)m).$$
\el \bpf This follows easily from  \eqref{1kr} and \eqref{kxt}.
\epf \noi
In addition we need a version of the Lemma for column addition. We leave it to the reader to fill in the details using  \eqref{2kr} and \eqref{kzt}.

 \subsection{$[X\ti \Z]$ as a graded poset}\label{gax}
We can regard $[X\ti \Z]$ as a  poset where $[{\mu},\ell] \le [{\gl},k]$ means that $[{\gl},k]$ can be obtained from $[{\mu},\ell]$ by applying a sequence of positive morphisms.  
We give the pair $({\gl},k)$ a degree by setting
\be \label{e67n}\deg({\gl},k)=|\gl| + k\in \Z.\ee  It is easily seen that this degree passes to equivalence classes. Thus 
$[X\ti \Z]= \bcu_{d\in \Z} [X\ti \Z]_{d}$ is a graded poset where 
$[X\ti \Z]_d = \{[{\gl},k]|\deg[{\gl},k] =d\}.$  
\bl  \label{wxt} For any $[\gl,  k] \in [X\ti \Z]$  we have $[\gl,  k]\le [\gl,  k+m]$
and $[\gl,  k]\le [\gl,  k+n].$
\el \bpf
If 
$\gl=(\gl_1, \gl_2, \ldots, \gl_n)$, define the partition $\mu$ by $\mu_1 = m$, 
$\mu_i = \gl_i$ for $i>1.$   Since $\mu$ is obtained from $\gl$ by adding boxes to the lowest row, we obtain 
$[\gl,  k]\le [\mu,  k] =[\mu^{-_\ttr} ,  k+m].$  By adding $\gl_n$ to the last part of  $\mu^{-_\ttr}$
and $\gl_i - \gl_{i+1}$ to the $i^{th}$ part for $i<n$, we see that 
$[\gl,  k]\le [\gl,  k+m]$.  The other statement is similar.
\epf \noi
The poset $ [X\ti Z]$ has an archimedian property.
\bc  \label{xxt}  For any $[\gl, in+jm]\in  [X\ti Z],$
we can find $K\in \Z$ such that $[\gl, in+jm]\le [\emptyset, Kmn]$.
\ec
\bpf Choose integers $k_0, k_1$ such that $i\le k_0m$ and $j\le k_1n$ and set $K = k_0+k_1 +1$.  Since 
$$ in +jm + (k_0m-i)n +(k_1n-j)n= (k_0+k_1)mn,$$ we have 
$$ [\gl, in+jm]\le [(m^n), in+jm]\le [\emptyset, Kmn].$$
\epf\noi 
\section{The Sergeev-Veselov Functor }\label{affc}

\subsection{The action on $\ttk^{n|m}$}\label{aog11}
In  \cite{SV101} Sergeev and Veselov defined 
an action of the Weyl groupoid  $\mathfrak{W}$  on $\ttk^{n|m}$ 
depending on a non-zero parameter $\gk$. 
We say $\gk$ is {\it special} if $\gk = \pm p/q$, where $p\in [m]$  and $q\in [n]$   with $p,q$ coprime.
If  $\gk = -1$ all orbits are described in \cite{M22}  Theorem 5.7. 
By  \cite{SV101} Theorem 3.2,  $\cW$ has an infinite orbit  if (and only if) $\gk$ is negative special. 
The proof of this statement reduces to the case $\gk=-m/n$.  We consider only the case
$\gk=-m/n$ and  
our notation  is slightly different from \cite{SV101}.  This will allow us to consider vectors and matrices with integer entries. The connection with the original notation is explained in Subsection \ref{cip}.  Starting in Subsection \ref{cpc} we assume that $m,n$ are coprime. The {\it Sergeev-Veselov Functor} 
$SV:\mathfrak  T_{iso}\lra \mathfrak \cS(\Z^{n|m})$
has already been defined in Subsection \ref{itr6}. 
For $\Gl$ as in \eqref{ezr}, set
\be \label{e73}
^L\Gl = {(a_1,\ldots, a_{n})}^\trans  , \quad ^R\Gl = ( b_1,\ldots ,b_m)
.\nn\ee 
We consider the $\mathfrak  T_{iso}$-orbit 
${\bf O}$ 
of 
   \be \label{e7}\Gl_0 = (m(n-1), \ldots,m,0| 0,n,\ldots,n(m-1))
= \sum_{i=1}^n  m(n-i)
\gep_i - \sum_{j=1}^m n(j-1)  \gd_j. 
\ee under the action from \eqref{e11}.

  \subsection{The augmented matrix $\widehat\A(\Gl)$}\label{fi}  
Let $\A(\Gl)$ be the matrix with entry in row $\gep_i$ and column $\gd_j$ given by 
\be \label{Adef}
\A(\Gl)_{i,j}= (\Gl,\gep_i-\gd_j) =a_i -b_j.
\ee 
For $a\in\ttk$, set 
$\bA=(a,\ldots, a|a,\ldots, a)$.  Suppose $\Gl$ is as in \eqref{ezr} and 
\be \label{eyr} \Gl' = 
(a'_1, \ldots ,a'_n|b'_1,\ldots,b'_{m}) .\ee

\bl \label{l73}
We have $\A(\Gl)=\A(\Gl') $ iff   $\Gl=\Gl' +\bA$ for some $a\in \ttk$.\el 
\bpf  If $\A(\Gl)=\A(\Gl') $, then $a= a_i -a'_i  =  b_j -b'_j$ is independent of $i,j$. 
\epf
\noi Since some information is lost in passing from $\Gl$ to $\A(\Gl)$,  we also use the (doubly) augmented matrix. 

\[\widehat \A(\Gl) =\begin{tabular}{|c||ccc|} \hline
&  &$^R\Gl$&
\\ \hline\hline &     &                          & \\ $^L\Gl$
&&$\A(\Gl)$
 & \\
 & &                              & \\
  \hline
\end{tabular}\]

\bl \label{fj} If $\Gl$  has integer entries and the entries  in  
$^L\Gl$ $($resp. $^R\Gl$$)$ form a complete set of representatives mod $n$ $($resp. mod $m$$)$, then the same is true for the entire $\cW$-orbit of $\Gl$.
\el
\bpf The condition is invariant under permutations of the entries in $\Gl$, so the result holds by \eqref{e9}.
\epf 
\bc \label{coD} If $\Gl\in {\bf O}$, then   $\A(\Gl)$ has at most one zero entry in each row and at most one zero entry in each column.
\ec
\noi 
In Example \ref{E2}, when $(n,m) =(2,3)$, we give several matrices $\widehat \A(\Gl)$. %
\subsection{Construction of $\Gl=x(\gl) \in \Z^{n|m}$ from $\Gl_0$ and $\gl \subseteq$ {\bf R}} \label{coL}
We use 
$\Gl_0$ as  in \eqref{e7}, and $\gl \subseteq$ {\bf R}, to construct certain elements 
$x(\gl)\in {\bf O}$.  However we point out that instead of $\Gl_0$, we could have used  
   \be \label{e7a}\Gl_0^{(k)} = \Gl_0 + k \bf{m}
\ee 
as a starting point for any $k\in \Z$.  Passing from $\Gl_0$ to $\Gl_0^{(k)} $,  the entries of $\Gl_0$ all increase by the same integer, so the difference between entries is unchanged. Thus $\A(\Gl_0^{(k)}) = \A(\Gl_0)$ and the key Equations \eqref{e67d} and \eqref{e67c} still hold when 
$ \Gl_0 $ is replaced by $\Gl_0^{(k)}$. In addition $\Gl_0^{(kn)}\in {\bf O}$   by Corollary \ref{aed}. 
 \\ \\
Given $\gl \subseteq$ {\bf R}, consider
\be \label{cd}  
\gl'=(\gl'_1,\ldots, \gl'_m) 
\mbox{ and } \lst =(\gl_n,\ldots, \gl_1)^\trans\ee 
as vectors.  Reversing the entries in $\gl$ means that when written as a column vector alongside $\gl \subseteq$ {\bf R}, the $i^{th}$ entry of $\lst$ gives the number of entries, $\gl_{n+1-i}$
in row $\gep_i $ of $\gl$. Define $\Gl =x(\gl)\in \Z^{n|m}$ by 
\be \label{e69}^L\Gl =^L\Gl_0 + n\lst \mbox{ and }^R\Gl={^R\Gl_0} + m\gl'.\ee
If   $\Gl_0$ and $\Gl$ are as in \eqref{e7} and \eqref{ezr} respectively, it follows that 
\be \label{e67a}
a_i=m(n-i)+ n\gl_{n+1-i}\ee
 and
\be \label{e67b}b_j=n(j-1) + m\gl'_j.\ee
Thus
\be \label{e67d}
a_i - 
a_k=m(k-i)+ n(\gl_{n+1-i}-\gl_{n+1-k})
\ee
and 
\be \label{e67c}b_{j}-b_{k} =n(j-k) + m(\gl'_{j}-\gl'_{k}).\ee
\bl \label{abc} If $\Gl=x(\gl)$ is as in \eqref{ezr} and $\ga = \gep_i-\gd_j$ is a positive root, 
\bi
\itema If $\gl\in X_\ga$ then $\Gl\in \Pi_\ga$.
\itemb If $\gl\in X_{-\ga}$ then $\Gl\in \Pi_{-\ga}$.
\itemc If 
$\gl \in X_{\pm\ga}$, then 
\be \label{abc1}x(t_{\pm\ga}(\gl))=\gt_{\pm\ga}(\Gl).\ee
\ei
\el
\bpf If $\gl\in X_{\ga}$, then $\gl_{n+1-i} = j-1$ and $\gl'_{j} = n-i$  
by Lemma \ref{rjx}.  Hence $a_i=b_j$ by Equations \eqref{e67a}  and \eqref{e67b}, that is $x(\gl)\in \Pi_\ga$ by  \eqref{e99}.  The proof of (b) is similar.
Finally, if 
$\gl \in X_{-\ga}$, and 
$\mu = t_{-\ga}(\gl)$, then by \eqref{qef}
 \[\mu_k= \left\{ \begin{array}
  {rcl}
  \gl_k -1  & \mbox{if} & k = {n+1-i} \\
\gl_k & \mbox{otherwise} &
\end{array} \right. \]
Similarly
 \[\mu'_k= \left\{ \begin{array}
  {rcl}
  \gl'_k -1  & \mbox{if} & k = j \\
\gl'_k & \mbox{otherwise} &
\end{array} \right. \] It follows from \eqref{e67a}  and \eqref{e67b}
that $x(\mu) = \gt_{-\ga}(\Gl).$ The proof of the other statement is similar.
\epf

\bc \label{aed} 
\bi \itema
If
$\gl \subseteq$ {\bf R}, then $x(\gl)\in   
{\bf O}$. 
\itemb With the notation of \eqref{e7a}, 
$\Gl_0^{(kn)}\in {\bf O}$ for all $k\in \Z$.  
\ei
\ec
\bpf For  (a) we use induction on $|\gl|$. The statement is trivial if $\gl = \emptyset$.  If 
$\mu =t_{\ga}(\gl)$
and  $x(\gl)\in   {\bf O},$  we 
obtain  $x(\mu)\in   {\bf O}$ by  \eqref{abc1}.   For (b), if  $\gl=(m^n)$ then by \eqref{e67a}
 and
\eqref{e67b}, $x(\gl)= \Gl_0  +\bf{mn}$, and this gives the result if $k=1$.    Now take any sequence of morphisms in $\mathfrak  T_{iso}$ that move 
$\Gl_0\in {\bf O}$ to 
$\Gl_0^{(n)}\in {\bf O}$ and apply their inverses in the reverse order to $\Gl_0$ reach 
$\Gl_0^{(-n)}$.  For $|k|>1$ the result follows by induction.  
\epf
\subsection{The coprime case} \label{cpc}
\subsubsection{ Converse to  Lemma \ref{abc} and exceptional cases} \label{spc}
We assume $m,n$ are coprime and obtain partial converses to the statements in Lemma \ref{abc}.  
 Suppose $\ga=\gep_i-\gd_j.$ Exceptional cases arise when ${ \boldsymbol \gl}\subseteq \gl $  or $\gl\in X^{\str}$.
\bl \label{c1} If 
$\Gl=x(\gl)\in \Pi_\ga$, then $\gl\in X_\ga$ unless 
${ \boldsymbol \gl}\subseteq \gl $ and $\ga=\gep_n-\gd_1$. 
\el
\bpf If $\Gl\in \Pi_\ga$, then 
$a_i=b_j$.  
Thus by 
\eqref{e67a}  and \eqref{e67b}
\be \label{gzr} n(\gl_{n+1-i}+1-j)= m(\gl'_j+i-n).\ee 
Since $i\in [m]$, $j\in [n]$, $0\le \gl'_j $ and $\gl_{n+1-i} \le m$, we have  
$$  \gl_{n+1-i}+1-j\le \gl_{n+1-i} \le m \mbox{ and } -n\le \gl'_j-n <  \gl'_j+i-n.$$
Thus, because $m,n$ are coprime  either both sides of \eqref{gzr} are zero, or both equal $mn$.   In the former case  we have 
\be \label{gjr}\gl_{n+1-i} = j-1 \quad\mbox{ and } \quad\gl'_{j} = n-i.\ee
  Thus  
$\gl\in X_\ga$ by Lemma \ref{rjx}. 
In the latter case  $i=n,$ $j=1$, $\gl_1=m$ and $\gl'_1=n$. 
This means ${ \boldsymbol \gl}\subseteq \gl$ and $\ga=\gep_n-\gd_1$. 
\epf
\bl \label{c2} If  $\Gl\in \Pi_{-\ga}$, then $\gl\in X_{-\ga}$ unless 
$\ga = \gep_1-\gd_{m}$ and $\gl_{n}= \gl'_m=0.$ 
\el
\bpf  
If $\Gl \in \Pi_{-\ga}$ then $a_i-b_j = (n-m)$, which 
by Equations \eqref{e67a}  and \eqref{e67b}
is equivalent to \be \label{gsr} n(\gl_{n+1-i}-j)= m(\gl'_j+i-n-1).\ee 
Both sides of \eqref{gsr} equal 0 or $-mn$.  In the first case 
we have 
$\gl_{n+1-i} = j$  and $\gl'_{j} =  {n+1-i},$ so  
 $\gl\in X_{-\ga}$  by Lemma \ref{rjs}.
In the second,  we deduce $i=1, j=m$ and $\gl_{n}= \gl'_m=0.$\epf
\subsubsection{ Comparison of $\A(\Gl)$ and $\A(\Gl^*)$} \label{spd}
\noi   Until the end of Subsection \ref{cpc}, we assume
\be \label{gjd}
\ga = \gep_i-\gd_{j}, \quad \gl\in X_{\ga}, \quad \gl^*=  t_{\ga}(\gl), \quad \Gl=x(\gl)\quad \mbox{ and } \quad\Gl^*= \gt_{\ga}(\Gl).\ee  
Thus $\Gl^* = x(\gl^*)$ by \eqref{abc1}. Note that $\Gl^*=\Gl + n\gep_i -m \gd_j$ by \eqref{e9}.  Also  $a_i-b_j =0$, see \eqref{ezr} and \eqref{exy}.  
Equation \eqref{gjr} holds by Lemma \ref{rjx},  
In Lemma \ref{obc} (resp. \ref{pbc}) we consider the possibility that $(\Gl^*,\gb) =0$ where $\gb $ is in the same row (resp. column) of {\bf R } as $\ga$.  
 Consider the roots $\ga^\ua, \ga^\ur$ which label the boxes immediately above and to the right of box $ \ga$.  These are given by 
$$\ga^\ua =\gep_{i-1}-\gd_{j}, \quad  \ga^\ur=\gep_i-\gd_{j+1}.$$
Of course $\ga^\ua, \ga^\ur$ are only defined if $i>1$ and $j<m$ respectively.  These conditions  will be assumed whenever the notation is used.

\bl \label{obc} If $(\Gl^*,\gb) =0$, for $\gb = \gep_i-\gd_k$ and $k\neq j$, then one of the following holds, 
\bi
\itema  $k=j+1$, $\gl'_{j}=\gl'_{j+1}$ and $\gb=\ga^\ur $.
\itemb $\ga = \gep_n-\gd_m$, $\gb = \gep_n-\gd_1$, 
$\gl'_{1}=n$, $\gl'_{m}=0$ 
and ${ \boldsymbol \gl}\subseteq \gl^*$. 
\ei
\el
\bpf  Since 
$(\Gl,\ga)=0,$ we have $(\Gl,\gb) =(\Gl,\gd_{j}-\gd_{k}) =b_{j}-b_{k} =n(j-k) + m(\gl'_{j}-\gl'_{k})$ by \eqref{e67c}. Thus as $(\Gl^*,\gb)=(\Gl,\gb) + n$, we have $(\Gl^*,\gb) =0$ iff
\be \label{hjd}
n(j-k+1) = m(\gl'_{k}-\gl'_{j}).\ee
First we show that $k \le j+1$.  Otherwise \eqref{hjd} implies that $k\ge m+j+1$, which is impossible since $j,k\in [m].$ 
If $k=j+1$, then \eqref{hjd} is zero and  we obtain the conclusion in (a). Finally if  $k \le j$, then $j+1=k+m$, which implies that $j=m, k=1, \gl'_{1}=n$  and $\gl'_{m}=0$. Since 
$\gl\in X_{\ga}$, we deduce that $i=n$. 
Since 
$\gl'_{1}=n$,
and $\ga = \gep_n-\gd_m$, $\gl^* = t_{\ga}(\gl)$ contains 
${ \boldsymbol \gl}$.
\epf

\bl \label{pbc} If $(\Gl^*,\gb) =0$, for $\gb = \gep_k-\gd_j$ and $k\neq i$, then one of the following holds, 
\bi
\itema  $k=i-1$, $\gl_{n+1-i}= \gl_{n+2-i}$ 
and $\gb=\ga^\ua $.
\itemb $\ga = \gep_1-\gd_1$, $\gb = \gep_n-\gd_1$, $\gl_{1}=m$, $\gl_{n}=0$ 
 and ${ \boldsymbol \gl}\subseteq \gl^*$. 
\ei
\el
\bpf  Since 
$(\Gl,\ga)=0,$ we have 
$(\Gl,\gb) =(\Gl,\gep_{k}-\gep_{i}) =a_k - a_i=m(i-k)+ n(\gl_{n+1-k}-\gl_{n+1-i})$ 
by \eqref{e67d}.  Thus as $(\Gl^*,\gb)=(\Gl,\gb) - m$, we have
$(\Gl^*,\gb) =0$ iff
\be \label{hjk}
m(k+1-i) = n(\gl_{n+1-k}-\gl_{n+1-i}).\ee
If $k+1<i $ we obtain 
$i\ge k+1+n$, which is impossible since $i,k\in [n].$ 
If $i=k+1$, then \eqref{hjk} is zero and  we obtain the conclusion in (a). Finally if  $k+1>i $, then $k+1=i+n$, which implies that $i=1, k=n, \gl_{1}=m$  and $\gl_{n}=0$. Since 
$\gl\in X_{\ga}$, we deduce that $j=1$. The conclusion in (b) follows. 
\epf

\subsubsection{A left inverse to $x$
}\label{daa}  
Let $Z$ be a set of zero entries in $\A(\Gl)$. We say a  descending path $p$ {\it supports } the entries in $Z$ if 
\bi\itema All the entries in $Z$ are above, to the right and adjacent to  $p$. 
\itemb Any other path with the property in (a) lies below and to the left of $p$.
\ei
There is always a unique path which supports more zero entries in $Z$ than any other.  
The entries under the path determine a partition $a(\Gl)$.  The corresponding set of zeroes is denoted  $Z(\Gl)$.  
\\ \\
Assume \eqref{gjd} holds and set $Z(\Gl)^\ga=Z(\Gl)\bsk \{\ga\}$.

\bl \label{le1} The set $Z(\Gl^*)$ satisfies \be \label{fjd}
Z(\Gl)^\ga \subseteq  Z(\Gl^*)\subseteq Z(\Gl)^\ga\cup \{\ga^\ua, \ga^\ur\}.\ee  
and is determined by \eqref{fjd} and the conditions
\bi \itema $\ga^\ur \in Z(\Gl^*)$ iff $\gl'_{j+1}= n-i.$  \itemb $\ga^\ua \in Z(\Gl^*)$ iff $\gl_{n+2-i}= j-1$. 
\ei 
\el

\bpf  
If $k\neq i$ and $\ell\neq j$, then $\A(\Gl)_{k,\ell}= \A(\Gl^*)_{k,\ell}$.  Also since $\A(\Gl)_{i,j} =0$, no other entry in  row $i$ or column $j$ of $\A(\Gl)$ can be zero  by Corollary  \ref{coD}.  This implies 
$Z(\Gl)^\ga \subseteq  Z(\Gl^*).$  
It remains to consider roots $\gb\neq \ga$ in the same row or column as $\ga$ such that 
$(\Gl^*,\gb)=0.$ We use \eqref{gjr}. If $\gb = \gep_i-\gd_k$, then by
Lemma \ref{obc}, either $k=j+1$, $n-i=\gl'_{j}=\gl'_{j+1}$ and $\gb=\ga^\ur $
or $\gb = \gep_n-\gd_1$. 
However in the latter case, the zero of $\A(\Gl^*)$ in box $\gb$ lies under the path defining 
$\gl^*$, so $\gb \notin Z(\Gl^*)$.  In the terminology of Subsection \ref{itr2} $\gb$ is an outer pseudo-corner of 
$\gl^*$. 
\\ \\
If $\gb = \gep_k-\gd_j$, then by
Lemma \ref{pbc} either $k=i-1$, $\gl_{n+1-i}= \gl_{n+2-i}$ 
and $\gb=\ga^\ua $ or 
$\ga = \gep_1-\gd_1$, $\gb = \gep_n-\gd_1$,
$\gl_{1}=m$, $\gl_{n}=0$ 
 and ${ \boldsymbol \gl}\subseteq \gl^*$. In the latter case, we  can show as above that $\gb \notin Z(\Gl^*)$.  
\epf

\noi Set 
$z(\gl)= \{\gb\in \Gd^+_1| \gl\in X_\gb\}$ and $z(\Gl)^\ga=z(\Gl)\bsk \{\ga\}$.

\bl \label{le7} The set $z(\gl^*)$ satisfies \be \label{jjd}
z(\gl)^\ga \subseteq  z(\gl^*)\subseteq z(\gl)^\ga\cup \{\ga^\ua, \ga^\ur\}.\ee  
and is determined by \eqref{jjd} and the conditions
\bi \itema $\ga^\ur \in z(\gl^*)$ iff $\gl'_{j+1}= n-i.$  \itemb $\ga^\ua \in z(\gl^*)$ iff $\gl_{n+2-i}= j-1$. 
\ei 
\el
\bpf
The upper boundary of a Young diagram in {\bf R} determines a path.  When $\ga$ is an outer corner  of $\gl$ we consider four consecutive steps in the paths for $\gl$  and $\gl^*$. For example if $\gl'_{j+1}= n-i$ and $\gl_{n+2-i}= j-1$
 we have the paths below with outer corners as shown. 
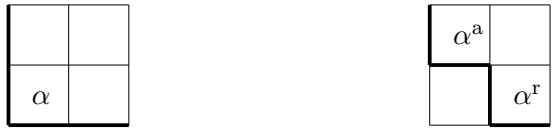
\Bc
\setlength{\unitlength}{0.8cm}
\begin{picture}(10,4)(-1,-.50)
\thinlines
  \linethickness{0.05mm}
\put(9.39,0.35){$\ga^\ur$}
\put(8.39,1.35){$\ga^\ua$}
\multiput(8.0,0)(1,0){3}{\line(0,1){2.0}}
 \multiput(8,0)(0,1){3}{\line(1,0){2.0}}
  \linethickness{0.5mm}
\put(1.0,0){\line(1,0){2}}
\put(1.0,2){\line(0,-1){2}}
  \linethickness{0.05mm}
\put(1.39,0.35){$\ga$}
\multiput(1.0,0)(1,0){3}{\line(0,1){2.0}}
 \multiput(1,0)(0,1){3}{\line(1,0){2.0}}
  \linethickness{0.5mm}
\put(8.0,1){\line(1,0){1}}
\put(9.0,0){\line(1,0){1}}
\put(8.0,2){\line(0,-1){1}}
\put(9.0,1){\line(0,-1){1}}
\end{picture}
\Ec
\noi 
The portion of the path on the left (resp. right) will be denoted d \underline{d} \underline{r}  r (resp. d  \underline{r}   \underline{d} r).  The underlined terms  are edges of $ \fB(\ga)$. 
Passing from the path for $\gl$ to the path in $\gl^*$, the order of the underlined terms is reversed. 
There are four possibilities and they determine $z(\gl^*)$ as in the table below. 

\[ \begin{tabular}{|c|c|c|c|c|} \hline
Path in $\gl$
 & r \underline{d} \underline{r}   d& r  \underline{d} \underline{r} r  & d \underline{d} \underline{r} d & d \underline{d} \underline{r} r \\ 
\hline Path in $\gl^*$
& r   \underline{r}   \underline{d} d & r \underline{r} \underline{d} r& d \underline{r}   \underline{d} d& d \underline{r}   \underline{d} r\\
  \hline $z(\gl^*)\bsk z(\gl)^\ga$
&$\emptyset$
&$\{\ga^\ur\}$& $\{\ga^\ua\}$& $\{\ga^\ur, \ga^\ua\}$\\
    \hline
\end{tabular}\]
The result follows from this.
\epf
\noi 
\bl \label{Th1} If  $a(\Gl)=\gl,$ then $\gl^*=a(\Gl^*)$. 

\el
\bpf  
The definition of 
$\gl=a(\Gl)$ depends on the set $Z(\Gl)$ which satisfies the recurrence in Lemma \ref{le1}.
On the other hand, the set $z(\gl)$ from  Lemma \ref{le7} satisfies the same recurrence.
\epf

\bc \label{ch1} If $\gl \subseteq$ {\bf R} and 
$\Gl=x(\gl)$, then $a(\Gl)=\gl.$  Thus $ax$ is the identity map on $X$.
\ec
\bpf  
Clearly the result holds if $|\gl|=0$.  So the result follows by induction on $|\gl|$ 
and Lemma \ref{Th1}. \epf

\subsection{Comparison of 
$x({\gl}), x({\gl}^\mr)$ and $x({\gl}^\omc)$} \label{rtn}
\noi If ${ \gl}\in X$ and 
 $\Gl= x(\gl)$, then the outer corners of  ${\gl}  $ occur in lower rows (moving from left to right). By Lemmas \ref{abc} and \ref{c1} the same is true for the zero entries in  $\A( \Gl)$ unless ${ \boldsymbol \gl}\subseteq \gl $ and $\ga=\gep_n-\gd_1$. 
However in the latter case, there is a unique descending path supporting all the zeroes in $\A(\nu(\Gl))$ (resp. $\A(\eta(\Gl))$). The partition deternined by this path is 
 ${\gl}^\mr $  (resp. ${\gl}^\omc $  ).  In this Subsection we note a simple relation between $x({\gl}), x({\gl}^\mr)$ 
and $x({\gl}^\omc)$. 
  See also Example \ref{E2}.
\\ \\
Suppose $\gl\in X(\rda)$.  
If $a_i(\gl)$, $b_j(\gl)$ are the expressions in \eqref{e67a} and 
\eqref{e67b}, then $a_i({\gl}^\mr)$, $b_j({\gl}^\mr)$ are given by
\be \label{e67f}
a_i({\gl}^\mr)=
m(n-i)+ n{\gl}^\mr_{n+1-i} = a_{\nu(i)}(\gl)-m \mbox{ for } i\in [n], i\neq 1.\ee
and 
\be \label{e67g}b_j({\gl}^\mr)=b_j(\gl) - m \mbox{ for } j\in [m].\nn\ee 
If $i=1$, then since ${\gl}^\mr_{\nu(1)} = {\gl}^\mr_n =0,$ \eqref{e67f} holds in this case also.  
Hence
\be \label{e67h}
x({\gl}^\mr)=
\sum_{i=1}^n  a_i({\gl}^\mr)\gep_i - \sum_{j=1}^m 
b_j({\gl}^\mr)\gd_j=-{\bf m} +\sum_{i=1}^n  
a_{\nu(i)}(\gl)\gep_i- \sum_{j=1}^m  b_j(\gl)\gd_j.\ee
Now \eqref{e67h} can be written in the form

\be \label{e67k}
x({\gl}^\mr)+{\bf m} =\nu(x({\gl})).\ee
Similarly we  have

\be \label{e67z}
x({\gl}^\omc)+{\bf n} =\eta(x({\gl})).\ee

\subsection{The  map 
$\hat x:[X\ti \Z] \lra\Z^{n|m}$} \label{eof}

We extend the map $x:X \lra \Z^{n|m}$ from Subsection \ref{coL}.  
Define a map 
$x:X\ti \Z \lra\Z^{n|m}$ by 
\be \label{xdef}
x(\gl,i{ n}+j{ m}) = 
{ \eta^{-i} \nu^{-j}}x(\gl) +i{\bf n}+j{\bf m}.\ee
Since the latter map has two arguments and the map $X \lra \Z^{n|m}$ has only one, there should be no danger of confusion with this notation. 
If  $\gl\in X(\rda)$, then  
\by \label{e01}  x({\gl}^\mr,i{ n}+(j+1){ m})
&=&   { \eta^{-i} \nu^{-j-1}}
x({\gl}^\mr)  +i{\bf n}+(j+1){\bf m}\nn\\
&=& 
{ \eta^{-i} \nu^{-j}}x(\gl) +i{\bf n}+j{\bf m}.
\ey  using \eqref{e67k}.  Thus 
$x({\gl}^\mr,k+m)=x(\gl,k)$.  Similarly since 
$x({\gl}^\omc,k+n)=x(\gl,k)$ if $\gl\in X(\cda)$,  
we see that  $x$ induces a map on equivalence classes, 
$\hat x:[X\ti \Z] \lra\Z^{n|m}$, $\hat x[\gl,in+jm] =x(\gl,in+jm)$.
\noi
\\ \\
We consider analogs of Lemmas \ref{abc}, \ref{c1} and \ref{c2} for the map $\hat x$.  In this situation there are no exceptional cases.
\bl \label{rex}   If 
$\ga \in \Gd_1^+$ and $s = [\gl,in+jm]$, then $s\in \cX_{\pm\ga}$ iff $\hat x(s)\in \Pi_{\pm\ga}$.
\el  \noi
\bpf Set $\gb= { \eta^{i} \nu^{j}}\ga$. 
First note that $\hat x(s) \in \Pi_{\pm\ga}$ iff $x(\gl) \in \Pi_{\pm\gb}$, by \eqref{xdef} and $W$-invariance of $(\;,\;)$. If $s\in \cX_{\pm\ga}$ it follows from Lemma \ref{abc}  that $\hat x(s) \in \Pi_{\pm\ga}$.  We prove the converse for positive roots.  If $x(\gl) \in \Pi_{\gb},$ then by Lemma \ref{c1}, either $\gl \in X_{\gb}$
or $\gb =\gep_n-\gd_1$ 
and ${ \boldsymbol \gl} \subseteq  \gl$.  In the former case $s\in \cX_{\ga}$. In the latter, by Lemma \ref{rey}, 
$\nu\gb$ 
 is an outer corner of 
${\gl}^\mr$.   It follows  that $s = [\gl^\mr,in+(j+1)m]\in \cX_{\ga}$.  For negative roots we can argue similarly using Lemma \ref{c2}, 
\epf
\noi 
The next result is crucial.
\bl \label{rjg} 
Set  $\gb= { \eta^{i} \nu^{j}}\ga.$  
 If $\gl\in X_{\pm\gb},$ 
then $\hat x(\gl,in+jm)\in 
\Pi_{\pm\ga}$  and 
\be \label{gjt} \hat x(T_{\pm\ga}[\gl,in+jm])= \gt_{\pm\ga}\hat x[\gl,in+jm].\ee
\el
\bpf  
If $\gl\in X_{\pm\gb},$ then 
 $\hat x[\gl,in+jm]\in \Pi_{\pm\ga}$ by Lemma \ref{abc} and \eqref{xdef}.
Therefore 
\by \label{e41}  \gt_{\pm\ga}\hat x[\gl,in+jm]   
&=&  \gt_{\pm\ga}[{ \eta^{-i} \nu^{-j}}x(\gl) +i{\bf n}+j{\bf m}]\nn\\
&=&  \gt_{\pm\ga}{ \eta^{-i} \nu^{-j}}x(\gl) +i{\bf n}+j{\bf m}.
\ey  
On the other hand by \eqref{dyx}, \eqref{xdef}   and then \eqref{abc1} 
\by 
\hat x(T_{\pm\ga}[\gl,in+jm]) &=&
\hat x(t_{\pm\gb}(\gl),i{n}+j{m}).\nn\\
&=&
{ \eta^{-i} \nu^{-j}}x(t_{\pm\gb}(\gl))+i{\bf n}+j{\bf m}.
\nn \\
&=&
{ \eta^{-i}\nu^{-j}}\gt_{\pm\gb}x(\gl)+i{\bf n}+j{\bf m}.
\nn
\ey  
and  \eqref{bc1} with $\Gl= x(\gl)$ shows that this is equal to \eqref{e41}.  
\epf \noi  Fix $\gl\in X, k \in \Z$  and set  $\Gl= \hat x(\gl,kmn)$ and set 
\be \label{fbc} Y(\gl) = \{\ga \in \Gd_1^+| \gl\in X_\ga\}, \quad Y(\Gl) = \{\ga \in \Gd_1^+| \Gl\in \Pi_\ga\}.
\ee

\bl \label{dbc} 
\bi\itemo
\itema If  $\gl = (m^r, 0^{n-r})$ with  $0\le r <n$ then  $Y(\gl)=Y(\Gl) = \{\gep_{n-r} -\gd_1\}$. 
\itemb
If $\gl = (s^{n})$ with  $0<  s <m$, then  $Y(\gl)=Y(\Gl) = \{\gep_{n} -\gd_{s+1}\}$.
\itemc
 If $\gl = (m^n)$ then  
$Y(\gl) = \emptyset $ and  $Y(\Gl) = \{\gep_{n} -\gd_1\}$.
\itemd In all other cases   $|Y(\gl)| > 1$ and $|Y(\Gl)| > 1$.
\iteme If $|Y(\gl)| = 1$ then $m$ divides $|\gl|$ or $n$ divides $|\gl|$.
\ei
\el
 
\bpf  The statements about $Y(\gl)$  are straightforward.  
Since $\A(\Gl) =\A(\Gl -k \bf{mn })$, we can assume $k=0$. Then 
 $Y(\gl) \subseteq Y(\Gl)$ by Lemma \ref{abc}. The rest follows from Lemma \ref{c1}.
\epf

\bl \label{zbc}  If $\Gl=\Gl_0 +k \bf{mn }$ with $k\in \Z$, then $|\hat x^{-1}(\Gl)|=1$.  
\el
\bpf Since $x(\emptyset, kmn) = \Gl$, $\hat x^{-1}(\Gl)$ is non-empty.
Using  \eqref{xdef} suppose that
$$\Gl=x(\gl,i{ n}+j{ m}) = 
{ \eta^{i} \nu^{j}}x(\gl) +i{\bf n}+j{\bf m}.$$ 
By a consideration of degrees, see \eqref{e67n} 
\be \label{xzef}
|\gl| + i{ n}+j{ m}= 
kmn.\ee  We may assume that $0\le j < n.$ From the above $ |Y(\Gl)| = |Y(x(\gl))| =1$.  Thus by Lemma \ref{dbc} (e) and \eqref{xzef} either 
 $i\equiv 0 \mod m$  or $j=0$  
and there are only $m+n$ possibilities for the pair $(i,j)$.  But    
$x^{-1}(\Gl)\subset X \ti Z$ is a union of equivalence classes as in Corollary \ref{vgs} and every equivalence class has cardinality $m+n$.  This proves the result. 
\epf
\bl \label{rqg} Suppose $s \in \cX_{-\ga},$ and $u =
T_{-\ga}(s).$ If  $|\hat x \hat x^{-1}(s)|=1$  then $|\hat x \hat x^{-1}(u)|=1$. 
\el
\bpf If  $\hat x(u')= \hat x  (u)$, then  $u, u' \in \cX_{\ga},$ by Lemma \ref{rex}. So  
$u' =
T_{-\ga}(s'),$ for some $s' \in \cX_{-\ga}.$ Thus  by Lemma \ref{rjg}
$\gt_{-\ga}\hat x(s) = \gt_{-\ga}\hat x(s')$. Since $ \gt_{-\ga}$ is bijective, we obtain $s= s'$ by hypothesis.
Therefore $u = u'$.   
\epf \noi  
\bt \label{rjt}
The  map 
$\hat x:[X\ti \Z] \lra\Z^{n|m}$ from \eqref{xdef} is injective with image  
${\bf O}$. 
\et \bpf  
If $\Gl = \hat x(\gl,k) \in {\bf O} \cap \Pi_{\gb}$, then 
$\gt_\gb(\Gl) \in \Im \hat x$ by Lemma \ref{rjg}.  This shows
${\bf O}\subseteq  \Im \hat x$.  
Next set $$S = \{s \in [X\ti \Z]| \hat x(s) \in  {\bf  O}\}.$$
If $s\in S\cap [X\ti \Z](\gb)$, then by Lemma 
\ref{rjg}, $ T_{\gb}(s) \in S$.   Since acts 
$\mathfrak T_{iso}$ transitively on $[X\ti \Z]$, we have $S= [X\ti \Z]$.  Thus 
$\Im \hat x  \subseteq {\bf O}$.
\\ \\
 To show injectivity we have to show   $|\hat x \hat x^{-1}(u)|=1$ for all  $u \in \cX.$
By Corollary  \ref{xxt}  and Lemma \ref{zbc}  we know  $|\hat x \hat x^{-1}(s)|=1$  for some $s \ge u$. 
Thus the result follows from Lemma \ref{rqg} by induction on deg $u -\deg s.$ 
\epf \noi 
Since $SV(\gr_{\pm\ga}) =\gt_{\pm\ga}$ and  $\F(\gr_{\pm\ga}) =T_{\pm\ga}$, we can 
write \eqref{gjt} in the form 
\be \label{gqt} 
SV(\gr_{\pm\ga})\hat x[\gl,k] =
\hat x(\F(\gr_{\pm\ga})[\gl,k]).\nn
\ee
Thus  $\hat x:[X \ti \Z] \lra{\bf O}$ is a
$\mathfrak T_{iso}$-equivariant bijection as claimed in  
Theorem \ref{iir}.

\bexa \label{E2}
{\rm 
\noi Take $(n,m) =(2,3)$.  We only indicate the zero entries in $\A(\Gl)$. 
 We have $\Gl_0 =(3.0|0,2,4)$ and  the  matrix $\hat \A(\Gl_0)$ is shown on the left below. On the right is the matrix $\hat \A(\Gl_1)$  arising from the morphism corresponding to the unique zero entry in 
$\A(\Gl_0)$. 
\Bc
\setlength{\unitlength}{1cm}
\begin{picture}(10,4)(-1,-.50)
\thinlines
  \linethickness{0.05mm}
\put(1.45,0.35){$0$}
\put(1.45,2.35){$0$}
\put(2.45,2.35){$2$}
\put(3.45,2.35){$4$}
\put(0.45,.35){$0$}
\put(0.45,1.35){$3$}
\multiput(1,0)(-1,0){2}{\line(0,1){3.1}}
\multiput(1.1,0)(1,0){4}{\line(0,1){3.1}}
 \multiput(0,2.1)(0,1){2}{\line(1,0){4.1}}
 \multiput(0,0)(0,1){3}{\line(1,0){4.1}}
  \linethickness{0.5mm}
\put(1.1,0){\line(1,0){3}}
\put(1.1,2){\line(0,-1){2}}
\thinlines
  \linethickness{0.05mm}
\put(9.45,.35){$0$}
\put(8.45,1.35){$0$}
\put(8.45,2.35){$3$}
\put(9.45,2.35){$2$}
\put(10.45,2.35){$4$}
\put(7.45,.35){$2$}
\put(7.45,1.35){$3$}
\multiput(8,0)(-1,0){2}{\line(0,1){3.1}}
\multiput(8.1,0)(1,0){4}{\line(0,1){3.1}}
 \multiput(7,2.1)(0,1){2}{\line(1,0){4.1}}
 \multiput(7,0)(0,1){3}{\line(1,0){4.1}}
  \linethickness{0.5mm}
  \put(8.1,1){\line(1,0){1}}
\put(9.1,1){\line(0,-11){1}}
 \put(9.1,0){\line(1,0){2}}
\put(8.1,2){\line(0,-11){1}}

\end{picture}
\Ec 
Given $\gl \subseteq$ {\bf R} as in Subsection \ref{coL} 
we can use  \eqref{e69}  to find $\Gl=x(\gl) \in {\bf O}$ as in Corollary \ref{aed}.  If 
$\gl=(3,1)$, the vectors $\gl', \lst$ from \eqref{cd} are given by $\gl'=(2,1,1)$ and $ \lst^\trans=(1,3)$.  Therefore $\Gl =(5,6|6,5,7).$ 
The matrix $\hat \A(\Gl)$
is given 
on the left below.  
The analysis in Subsection \ref{daa} depended on the unique 
 descending path  that supports  the maximum number of zero entries in the augmented matrix.  Note for this matrix $\hat \A(\Gl)$, 
there is no descending path supporting {\it all } the zeroes, compare Lemma \ref{c1}. The descending path supporting the upper zero
corresponds to the partition $\gl$.  We have ${\gl}^\mr=(1,0)$.  
After rotating rows in the diagram on the left, we obtain the diagram on the right where the path corresponds to ${\gl}^\mr$.  Note that all entries in this last matrix are the same as those in $\hat \A(\Gl_1)$, except that   the entries in $^L\Gl_1$ and $^R\Gl_1$ have been increased by $m=3$.  In the notation of 
 \eqref{e67k}
\be \label{gl9}
x(1,0)+{\bf 3} =\nu(x(3,1)).\ee

\Bc
\setlength{\unitlength}{1cm}
\begin{picture}(10,4)(-1,-.50)
\thinlines
  \linethickness{0.05mm}
\put(1.45,0.35){$0$}
\put(2.45,1.35){$0$}
\put(1.45,2.35){$6$}
\put(2.45,2.35){$5$}
\put(3.45,2.35){$7$}
\put(0.45,.35){$6$}
\put(0.45,1.35){$5$}
\multiput(1,0)(-1,0){2}{\line(0,1){3.1}}
\multiput(1.1,0)(1,0){4}{\line(0,1){3.1}}
 \multiput(0,2.1)(0,1){2}{\line(1,0){4.1}}
 \multiput(0,0)(0,1){3}{\line(1,0){4.1}}
\thinlines

  \linethickness{0.05mm}
\put(8.45,1.35){$0$}
\put(9.45,.35){$0$}
\put(8.45,2.35){$6$}
\put(9.45,2.35){$5$}
\put(10.45,2.35){$7$}
\put(7.45,.35){$5$}
\put(7.45,1.35){$6$}
\multiput(8,0)(-1,0){2}{\line(0,1){3.1}}
\multiput(8.1,0)(1,0){4}{\line(0,1){3.1}}
 \multiput(7,2.1)(0,1){2}{\line(1,0){4.1}}
 \multiput(7,0)(0,1){3}{\line(1,0){4.1}}
  \linethickness{0.5mm}

  \put(8.1,1){\line(1,0){1}}
\put(9.1,1){\line(0,-11){1}}
 \put(9.1,0){\line(1,0){2}}
\put(8.1,2){\line(0,-11){1}}
  \put(1.1,2){\line(1,0){1}}
\put(2.1,2){\line(0,-11){1}}
 \put(2.1,1){\line(1,0){2}}
\put(4.1,1){\line(0,-11){1}}

\end{picture}
\Ec For an analog of \eqref{gl9} involving $\eta$ and where a column is deleted, note that, 
$$
x(2,0)+{\bf 2} =\eta(x(3,1)).$$

}\eexa

\noi



\section{A conjecture on infinite orbits}\label{cio}
The orbit representative $\Gl_0$  from 
 \eqref{e7} has some special properties that enable us to relate the orbit {\bf O} to Young diagrams.   
We conjecture that {\bf O} is 
 the unique infinite orbit up to some inherent symmetries.
\bco \label{t2} Up to a  translation by a suitable vector $\bA$ 
any infinite orbit  for the $\mathfrak{W}$ action on $\ttk^{n|m}$  has a representative $\Gl_0$  of the form 
 \eqref{e7}.
\eco \noi 
\subsection{Motivation}\label{cip}
First we explain how our notation related to 
 \cite{SV101}. 
We say $\gk$ is {\it special} if $\gk = \pm p/q$, where $p\in [m]$  and $q\in [n]$   with $p,q$ coprime.
 The following result is contained in \cite{SV101}  Theorem 1.1.  
\begin{thm} 
\label{finite}
All the orbits are finite if and only if $\kappa$ is not negative special.
\end{thm}
\noi Since we are interested in infinite orbit we assume $\gk = -p/q$ is rational. The action from  \cite{SV101}  Equation (9)  is given as follows.   We assume that $\ga=\gep_i-\gd_j$ and omit the functor.  Then on the hyperplane $x_i=y_j$ we have 
$$\gr_\ga(x_1, \ldots, x_n, y_1, \ldots, y_m) = (x_1, \ldots, x_i+1, \ldots, x_n, y_1, \ldots, y_j-\gk, \ldots, y_m).$$ We leave it to the reader to determine the inverse morphism.  
Introduce new coordintates $X_i =qx_i, Y_j = qy_j$.  Then the action is given by 
$$\gr_\ga(X_1, \ldots, X_n, Y_1, \ldots, Y_m) = (X_1, \ldots, X_i+q, \ldots, X_n, Y_1, \ldots, Y_j+p, \ldots, Y_m).$$ If $p=m$ and $q=n$ we have the action 
defined by the functor $SV$ from \eqref{e11}- \eqref{e9}.  In general let 
$v_\ga^{p,q} = q\gep_i -p\gd_j$
 and 
\be \label{e95}\Pi_{\ga}^{p,q} = \{\Gl\in \ttk^{n|m}| (\Gl, \ga)=0\},\quad \Pi_{-\ga}^{p,q} = \{\Gl\in \ttk^{n|m}| (\Gl,\ga)=q-p\}.\ee
Then  we have a functor $SV^{p,q}:\mathfrak  T_{iso}\lra \mathfrak \cS(\Z^{n|m})$ given by 
$SV^{p,q}(\pm \ga) =\Pi_{\pm\ga}^{p,q}$, and 
\be \label{e90}SV^{p,q}(\gr_{\pm\ga})(\Gl) = \Gl \pm v_\ga^{p,q}.\ee 
Now we change our point of view.   Fix $n'\le n$ and $m'\le m$ , with $n',m'$ coprime and let $\mathfrak  T_{iso}^{n',m'}$ be the subgroupoid of $\mathfrak  T_{iso}$ with base $\{\gep_i-\gd_j|i\in [n'], j\in [m']\}$ and restrict  the 
action $SV:\mathfrak  T_{iso}\lra \mathfrak \cS(\Z^{n|m})$ given by 
\eqref{e11} to $\mathfrak  T_{iso}^{n',m'}$.   Then from Theorem \ref{finite} we have  
\bl \label{Ade} If $n'< n$ or $m'< m$,  then all orbits of  $\mathfrak  T_{iso}^{n',m'}$ on $\Z^{n|m}$ are finite.
\el
\bpf  The parameter $-n/m$ is not negative special for the action of $\mathfrak  T_{iso}^{n',m'}$.
\epf
\noi 
The parameter governing the situation considered in Section \ref{affc} is thus on the borderline for the existence of an infinite orbit.  On the other hand it is also borderline for 
Conjecture \ref{t2}: if we consider the action of  $\mathfrak  T_{iso}^{n|m+1}$ on $\mathfrak \cS(\Z^{n|m+1})$ with parameter $-n/m$, then we obtain an infinite orbit using only the first $m$ columns (and all $n$ rows).  The entries in the last column can be arbitrary.   Similar remarks apply to the the action of  $\mathfrak  T_{iso}^{n+1|m}$ on $\mathfrak \cS(\Z^{n+1|m})$ with parameter $-n/m$


\begin{thebibliography}{99}


\bib{M101}{book}{
   author={Musson, Ian M.},
   title={Lie superalgebras and enveloping algebras},
   series={Graduate Studies in Mathematics},
   volume={131},
   publisher={American Mathematical Society},
   place={Providence, RI},
   date={2012},
   pages={xx+488},
   isbn={978-0-8218-6867-6},
   review={\MR{2906817}},
}

\bib{M22}{article}{
   author={Musson, Ian M.},
   title={On the geometry of some algebras related to the Weyl groupoid},
   conference={
      title={Recent advances in noncommutative algebra and geometry},
   },
   book={
      series={Contemp. Math.},
      volume={801},
      publisher={Amer. Math. Soc., [Providence], RI},
   },
   date={2024},
   pages={155--177},
   review={\MR{4756383}},
   doi={10.1090/conm/801/16039},
}
		


\bib{M23}{article}{
      author={Musson, Ian M.},  
       title={Young diagrams, Borel subalgebras and Cayley graphs, arXiv:2412.12141 
},
        type={In preparation
},
       date={2024},
}

\bib{S}{article}{
   author={Serganova, Vera},
   title={On generalizations of root systems},
   journal={Comm. Algebra},
   volume={24},
   date={1996},
   number={13},
   pages={4281--4299},
   issn={0092-7872},
   review={\MR{1414584}},
   doi={10.1080/00927879608825814},
}





\bib{SV1}{article}{
   author={Sergeev, A. N.},
   author={Veselov, A. P.},
   title={Deformed quantum Calogero-Moser problems and Lie superalgebras},
   journal={Comm. Math. Phys.},
   volume={245},
   date={2004},
   number={2},
   pages={249--278},
   issn={0010-3616},
   review={\MR{2039697}},
   doi={10.1007/s00220-003-1012-4},
}




\bib{SV2}{article}{
   author={Sergeev, A. N.},
   author={Veselov, A. P.},
   title={Grothendieck rings of basic classical Lie superalgebras},
   journal={Ann. of Math. (2)},
 volume={173},
   date={2011},
   number={2},
   pages={663--703},
   issn={0003-486X},
}
\bib{SV101}{article}{
   author={Sergeev, A. N.},
   author={Veselov, A. P.},
   title={Orbits and invariants of super Weyl groupoid},
   journal={Int. Math. Res. Not. IMRN},
   date={2017},
   number={20},
   pages={6149--6167},
   issn={1073-7928},
   review={\MR{3712194}},
   doi={10.1093/imrn/rnw182},
}
		





\bib{We}{article}{
   author={Weinstein, Alan},
   title={Groupoids: unifying internal and external symmetry. A tour through
   some examples},
   journal={Notices Amer. Math. Soc.},
   volume={43},
   date={1996},
   number={7},
   pages={744--752},
   issn={0002-9920},
   review={\MR{1394388}},
}
	
\end{thebibliography}
\end{document}